\documentclass[12pt]{article}
\usepackage{latexsym, amssymb,amsfonts,amsmath, psfrag,eepic,colordvi,epsfig,ifpdf}
\newtheorem{theo}{Theorem}[section]

\def\qed{\hfill \rule{4pt}{7pt}}
\def\pf{\noindent {\it Proof.} }

\begin{document}
\parskip 6pt

\pagenumbering{arabic}
\def\sof{\hfill\rule{2mm}{2mm}}
\def\ls{\leq}
\def\gs{\geq}
\def\SS{\mathcal S}
\def\qq{{\bold q}}
\def\MM{\mathcal M}
\def\TT{\mathcal T}
\def\EE{\mathcal E}
\def\lsp{\mbox{lsp}}
\def\rsp{\mbox{rsp}}
\def\pf{\noindent {\it Proof.} }
\def\mp{\mbox{pyramid}}
\def\mb{\mbox{block}}
\def\mc{\mbox{cross}}
\def\qed{\hfill \rule{4pt}{7pt}}
\def\block{\hfill \rule{5pt}{5pt}}
\begin{center}
{\Large {\bf Noncrossing Linked Partitions and \\[5pt]Large
$(3,2)$-Motzkin Paths}}\\

\vskip 6mm \vskip 6mm

{\small William Y.C. Chen$^1$, Carol J. Wang$^2$\\[2mm]
$^1$Center for Combinatorics, LPMC-TJKLC\\
Nankai University\\
Tianjin 300071, P.R. China \\[3mm]

$^2$Department of Applied Mathematics \\
Beijing Technology and
Business University\\ Beijing 100048, P.R. China \\[3mm]
$^1$chen@nankai.edu.cn, $^2$wang$\_$jian@th.btbu.edu.cn}

\vskip 6mm
\end{center}

\vskip10mm
\noindent {\bf Abstract.} Noncrossing linked partitions arise in the
study of certain transforms in free probability theory. We explore the
connection between noncrossing linked partitions and colored Motzkin paths.
A $(3,2)$-Motzkin path can be viewed as a colored
Motzkin path in the sense that there are three types of level steps and
two types of down steps. A large $(3,2)$-Motzkin path is defined to be
a $(3,2)$-Motzkin path for which  there are only two types of level steps on the
$x$-axis. We establish a one-to-one correspondence between the set of noncrossing linked
partitions of $[n+1]$ and the set of large $(3,2)$-Motzkin paths of
length $n$. In this setting, we get a
 simple explanation of the
well-known relation between the large and the little Schr\"oder numbers.



\vskip10mm

\section{Introduction}

The notion of noncrossing linked partitions was introduced by
Dykema \cite{Dyke07} in the study of free probabilities. He showed
that the generating function of the number of noncrossing linked
partitions of $[n+1] =\{1,2,\ldots,n+1\}$ equals
\begin{equation}\label{Fx}F(x)=\sum_{n=0}^{\infty}f_{n+1}x^n
=\frac{1-x-\sqrt{1-6x+x^2}}{2x}.\end{equation} This implies that the
number of noncrossing linked partitions of $[n+1]$ is equal to the number of large Schr\"{o}der paths of length $2n$, namely, the
$n$-th large Schr\"{o}der number $S_n$.
 A large Schr\"{o}der
path of length $2n$ is a lattice path from $(0,0)$ to $(2n,0)$
consisting of up steps $(1,1)$, level steps $(2,0)$ and down
steps $(1,-1)$ and never lying below the $x$-axis. The first few
values of  $S_n$'s are
$1,2,6,22,90,394,1806,\ldots$. The sequence of the large Schr\"{o}der numbers is
listed as  entry A006318 in OEIS \cite{EIS}. Chen, Wu and
Yan~\cite{CWY08} established a bijection between the set of
noncrossing linked partitions of $[n+1]$ and the set of large
Schr\"{o}der paths of length $2n$.

Motivated by the correspondence between noncrossing partitions and $2$-Motzkin paths, we
are led to the question whether there is any connection between
noncrossing linked partitions and colored Motzkin paths.
It is known that  little Schr\"oder paths of length $2n$ equals the
 number of $(3,2)$-Motzkin paths of length $n-1$. Recall that a little Schr\"oder
 path of length $2n$
 is a lattice path from $(0,0)$ to $(2n,0)$
  consisting of up steps $(1,1)$, down steps $(1,-1)$ and level steps $(2,0)$ not lying below the $x$-axis with the additional condition that there are no level steps on the $x$-axis. The number of such paths of length $2n$ is referred to as the little Schr\"oder number $s_n$. Yan \cite{HFYan07} found a bijective
 proof of this fact. Since the large Schr\"oder numbers and the little Schr\"oder numbers
 are related by a factor of two, we see that the number of noncrossing linked partitions
 of $[n+1]$ is numerically related to the number of $(3,2)$-Motzkin paths
 of length $n$.

 Indeed, the main result of this paper is to introduce a
 class of Motzkin paths, which we call the large $(3,2)$-Motzkin paths,  such that
  noncrossing linked partitions of $[n+1]$ are in one-to-one correspondence with large
 $(3,2)$-Motzkin paths of length $n$. By examining the connection between large $(3,2)$-Motzkin paths and
 the ordinary $(3,2)$-Motzkin paths, we immediately get the relation between the
 large and the little Schr\"oder numbers.


Let us recall some terminology. A $(3,2)$-Motzkin path of
length $n$ is a lattice path from $(0,0)$ to $(n,0)$ consisting of
up steps $u=(1,1)$, level steps $(1,0)$ and down steps $(1,-1)$ with
each down step receiving one of the two colors $d_1, d_2$, and each
level step receiving one of the three colors $l_1, l_2, l_3$. Let $m_n$ denote the
$n$-th $(3,2)$-Motzkin number, that is, the number of
$(3,2)$-Motzkin paths consisting of $n$ steps, or of length $n$. A
large $(3,2)$-Motzkin path is a $(3,2)$-Motzkin path for
which each level step at the $x$-axis receives only one of the two colors
$l_1$ or $l_2$. An elevated large $(3,2)$-Motzkin path is
defined as a large $(3,2)$-Motzkin path that does not touch the
$x$-axis except for the origin and the destination.
 Denote the set
of large $(3,2)$-Motzkin paths by $L$ and the set of large
$(3,2)$-Motzkin paths of length $n$ by $L(n)$.
Meanwhile, we use  $L_n$ to denote
the number of paths in $L(n)$.

It can be shown that the generating function $$L(x)=\sum_{n=0}^{\infty}L_nx^n$$ satisfies the functional equation
\begin{equation}\label{lx}
L(x)=1+2xL(x)+2x^2M(x)L(x),
\end{equation}
where
\begin{equation}\label{mx}
M(x)=\sum_{n=0}^{\infty}m_nx^n
=\frac{1-3x-\sqrt{1-6x+x^2}}{4x^2}
\end{equation}
 is
the generating function of the $(3,2)$-Motzkin numbers. From (\ref{lx}) and (\ref{mx}),
it follows that $L(x)=F(x)$.
This yields that
\begin{equation}\label{Lf2}
L_n=f_{n+1}.\end{equation}

Examining the connection between the large $(3,2)$-Motzkin
paths and
ordinary $(3,2)$-Motzkin paths, we are
led to a simple explanation of
the following
relation:
\begin{equation}\label{Lm2}
L_n=2m_{n-1}.
\end{equation}
In fact, the argument for the above relation is
essentially the same argument for the relation
(\ref{lx}) between $L(x)$ and $M(x)$.
Since the little Schr\"oder number is equal to the $(3,2)$-Motzkin number, see Chen, Li, Shapiro and
Yan \cite{CLSY07} and Yan \cite{HFYan07},  (\ref{Lm2}) is equivalent to the
well-known relation $S_n=2s_n$, which has been proved combinatorially by Shapiro and Sulanke \cite{ShSu00}, Deutsch \cite{Deut01}, Gu, Li and Mansour \cite{GLM08} and Huq \cite{Huq09}.

\section{Noncrossing Linked Partitions}

 We begin with an overview of
noncrossing linked partitions. In particular,
we shall give a description of the linear
representation of a linked partition.
A  linked partition $\pi=\{B_1,B_2,\ldots,B_k\}$ of $[n]$ is a collection
 of nonempty subsets, called blocks, of $[n]$,  such that
the union of $B_1, B_2, \ldots, B_k$ is $[n]$ and any two distinct blocks
of $\pi$ are nearly disjoint. Two blocks
 $B_i$ and $B_j$  are said to be
  nearly disjoint if for every $k \in B_i\cap B_j$, one of
the following conditions holds:\\[4pt]
(a) $k=\min(B_i), |B_i|>1$ and $k\neq \min(B_j)$, or\\[4pt]
(b) $k=\min(B_j), |B_j|>1$ and $k\neq \min(B_i)$.

We say that $\pi$ is a noncrossing linked partition if in addition, for any two distinct blocks $B_i,B_j\in \pi$, there does not exist $i_1,j_1\in B_i$ and $i_2,j_2\in B_j$ such that
$i_1<i_2<j_1<j_2$. Let $NCL(n)$ denote the set of noncrossing
linked partitions of $[n]$.

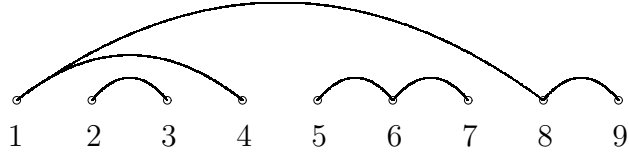
\begin{figure}[h]
\begin{center}
\setlength{\unitlength}{0.2cm}
\begin{picture}(40,10)
\multiput(0,3)(5, 0){9}{\circle{0.5}}\put(-.6,0){$1$}\put(4.6,0){$2$}
\put(9.6,0){$3$}\put(14.6,0){$4$}\put(19.6,0){$5$}\put(24.6,0){$6$}\put(29.6,0){$7$}
\put(34.6,0){$8$}\qbezier(0,3)(17.5,16)(35,3)\qbezier(0,3)(7.5,9)(15,3)
\qbezier(5,3)(7.5,6)(10,3)\qbezier(20,3)(22.5, 6)(25,3)
\qbezier(25,3)(27.5,6)(30,3)
\qbezier(35,3)(37.5,6)(40,3)\put(39.6,0){$9$}
\end{picture}
\caption{The linear representation of $\pi=\{1,4,8\}\{2,3\}\{5,6\}\{6, 7\}\{8,9\}$}.\label{NCLP}
\end{center}
\end{figure}

 In this paper, we shall adopt the linear representation  of noncrossing linked partitions, see Chen, Wu and
Yan \cite{CWY08}. For a noncrossing linked partition $\pi$ of $[n]$, first we draw $n$ vertices on a horizontal line with points or vertices $1,2,\ldots,n$ arranged in increasing order. For each block $B=\{i_1,i_2,\ldots,i_k\}$ with $i_1=\min(B)$ and $k\geq 2$, draw an arc between $i_1$ and any other vertex in $B$.  We shall use a pair $(i,j)$ to denote an arc between $i$ and $j$, where we assume that $i<j$. For example, the linear representation of
$\pi=\{1,4,8\}\{2,3\}\{5,6\}\{6, 7\}\{8,9\}\in NCL(9)$ is shown in Figure \ref{NCLP}.

Here is the main result of this paper.

\begin{theo}
There is a bijection between the set of large $(3,2)$-Motzkin paths
 of length $n$ and the set of noncrossing linked partition of $[n+1]$.
\end{theo}

\pf We describe a map $\varphi$ from $L(n)$ to $NCL(n+1)$ in terms of a recursive procedure.
Let $P\in L(n)$. We wish to construct a noncrossing linked partition $\pi= \varphi(P)$.

If $P=\emptyset$, then set $\varphi(P)=\{1\}$.
For $n\geq 1$, write $P=P_1P_2\cdots P_k,$ where each
$P_i$ is a nonempty elevated large $(3,2)$-Motzkin paths of length $p_i$.  We consider the following  cases.


{\setlength{\leftmargin}{2cm}
\begin{enumerate}


\item[Case 1.]
\begin{itemize}\item[(i)] If $P_i=l_1$, then set
$\varphi(P_i)=\{1,2\}$.
\item[(ii)] If $P_i=l_2$, then set $\varphi(P_i)=\{1\}\{2\}$.
\end{itemize}

\item[Case 2.] \begin{itemize}\item[(i)] If $P_i=uP_{c}d_1$, where $P_{c}\in L$, that is, $P_c$ is a large $(3,2)$-Motzkin path,
then set  $$\varphi(P_i)=\{1,p_i,p_i+1\}\cup \varphi(P_{c}),$$
see Figure~\ref{Case2a} for an illustration of this operation.

\begin{figure}[h]
\begin{center}
\setlength{\unitlength}{0.5cm}
\begin{picture}(18,2.5)
\put(0,0){\circle*{.1}}\put(0.26,0.6){$u$}\put(0,0){\line(1,1){1}}\put(1,1){\circle*{.1}}\qbezier[30](1,1)(3,1)(5,1)
\qbezier[30](1,1)(3,2)(5,1)\put(5,1){\circle*{.1}}\put(5,1){\line(1,-1){1}}\put(6,0){\circle*{.1}}
\put(2.1,2){$P_{c}\in L$}\put(5.5,0.6){$d_1$}\qbezier[40](-1,0)(3,0)(7,0)

\put(7.5,.5){\vector(1,0){1.5}}\put(8.1,0.7){$\varphi$}

\put(10,0.3){\circle{.15}}\put(16,0.3){\circle{.15}}\put(17.5,0.3){\circle{.15}}
\put(9.45,-.05){$1$}\put(15.85,-.35){$p_i$}\put(17,-.35){$p_i+1$}
\qbezier(10,.3)(13,2.3)(16,.3)
\put(10,-.35){\dashbox{.2}(5,.65){$\scriptstyle \varphi(P_{c})$}}
\qbezier(10,.3)(13.75,3.3)(17.5,.3)
\end{picture}
\caption{Case 2 (i).}\label{Case2a}
\end{center}
\end{figure}
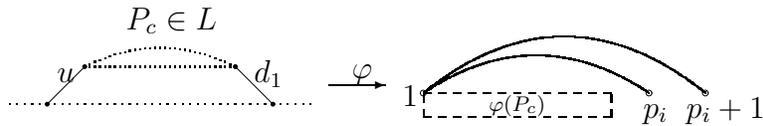

\item[(ii)] If $P_i=uP_cd_1$ and $P_c\in \overline{L}$, that is, $P_c$ is a $(3,2)$-Motzkin path  with at least one $l_3$ step on the $x$-axis. Then set $$P_c=P_c^{(1)}l_3P_c^{(2)}l_3\cdots
l_3P_c^{(k)}, \ k\geq 2, $$ where $P_c^{(i)}\in L$ is of length $t_i \geq 0$.
 We proceed to construct $\varphi(P_i)$ via the following steps. Let
$$\tau(P_c^{(i)})=\{1,t_i + 2\}\cup \varphi(P_c^{(i)}),\ i=1,2,\ldots, k.$$
For $i=1,2,\ldots, k-1$, merge the last vertex  $t_i + 2$ of $\tau(P_c^{(i)})$ and the first vertex $1$ of $\tau(P_c^{(i+1)})$ and relabel the vertices by $\{1,2,\ldots,p_i\}$ in increasing order. Denote the resulting noncrossing linked partition by $\omega(P_c)$. Then set
$$\varphi(P_i)=\{1,p_i+1\}\cup\omega(P_c).$$
An illustration of the above construction is given in Figure~\ref{Case2b}.

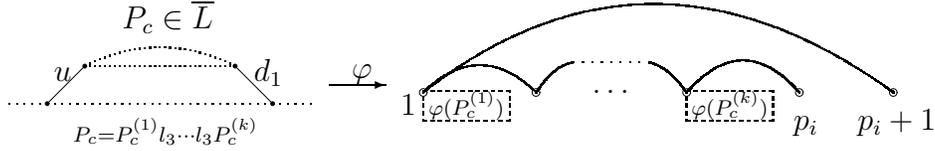
\begin{figure}[h]
\begin{center}
\setlength{\unitlength}{0.5cm}
\begin{picture}(23,4)
\put(0,1){\circle*{.1}}\put(0.18,1.6){$u$}\put(0,1){\line(1,1){1}}\put(1,2){\circle*{.1}}\qbezier[30](1,2)(3,2)(5,2)
\qbezier[30](1,2)(3,3)(5,2)\put(5,2){\circle*{.1}}\put(5,2){\line(1,-1){1}}\put(6,1){\circle*{.1}}
\put(0.7,0){$\scriptstyle P_c=P_c^{(1)}l_3\cdots
l_3P_c^{(k)}$}\put(2,3){$P_c\in \overline{L}$}
\put(5.5,1.6){$d_1$}\qbezier[40](-1,1)(3,1)(7,1)

\put(7.5,1.5){\vector(1,0){1.5}}\put(8.1,1.7){$\varphi$}

\multiput(10,1.3)(3,0){2}{\circle{.2}}
\multiput(17,1.3)(3,0){2}{\circle{.2}}\put(22.5,1.3){\circle{.2}}
\put(14.5,1.3){$\ldots$}

\qbezier(10,1.3)(11.5,2.75)(13,1.3)\qbezier(13,1.3)(13.5,2.05)(14,2.1)
\qbezier[8](14,2.1)(15,2.1)(16,2.1)\qbezier(16,2.1)(16.5,2)(17,1.3)
\qbezier(17,1.3)(18.5,3)(20,1.3)\qbezier(10,1.3)(16,6)(22.5,1.3)

\put(9.4,0.7){$1$}
\put(19.85,0.3){$p_i$}\put(21.55,0.3){$p_i+1$}

\put(10,0.55){\dashbox{.1}(2.3,.8){{$\scriptstyle
\varphi(P_c^{(1)})$}}}

\put(17,0.55){\dashbox{.1}(2.3,.8){{$\scriptstyle
\varphi(P_c^{(k)})$}}}

\end{picture}
\caption{Case 2 (ii).}\label{Case2b}
\end{center}
\end{figure}
\end{itemize}

\item[Case 3.] \begin{itemize}\item[(i)] If $P_i=uP_{c}d_2$, where $P_{c}\in L$,
then  set $$\varphi(P_i)=\left(\{1,p_{i}+1\}\{p_i\}\right)\cup
\varphi(P_{c}).$$
 Figure~\ref{Case3a} is an illustration of this operation.

\begin{figure}[h]
\begin{center}
\setlength{\unitlength}{0.5cm}
\begin{picture}(18,3)
\put(0,0){\circle*{.1}}\put(0.26,0.6){$u$}\put(0,0){\line(1,1){1}}\put(1,1){\circle*{.1}}\qbezier[30](1,1)(3,1)(5,1)
\qbezier[30](1,1)(3,2)(5,1)\put(5,1){\circle*{.1}}\put(5,1){\line(1,-1){1}}\put(6,0){\circle*{.1}}
\put(2.1,2){$P_{c}\in L$}\put(5.5,0.6){$d_2$}\qbezier[40](-1,0)(3,0)(7,0)

\put(7.5,.5){\vector(1,0){1.5}}\put(8.1,0.7){$\varphi$}

\put(10,0.3){\circle{.15}}\put(16,0.3){\circle{.15}}\put(17.5,0.3){\circle{.15}}
\put(9.45,-.05){$1$}\put(15.85,-.35){$p_i$}\put(17,-.35){$p_i+1$}

\put(10,-.35){\dashbox{.2}(5,.65){$\scriptstyle \varphi(P_{c})$}}
\qbezier(10,.3)(13.75,3.3)(17.5,.3)
\end{picture}
\caption{Case 3 (i).}\label{Case3a}
\end{center}
\end{figure}
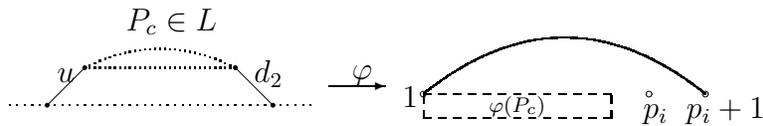

\item[(ii)] If $P_i=uP_cd_2$ and $P_c\in \overline{L}$, then set
$$P_c=P_c^{(1)}l_3P_c^{(2)}l_3\cdots l_3P_c^{(k)}, \ k\geq 2,$$
where $P_c^{(i)}\in L$ is of length $t_i \geq 0$. Denote by $\tau(P_c^{(i)})$
the noncrossing linked partition
$\{1,t_i+2\}\cup \varphi(P_c^{(i)}),\ i=2,3,\ldots, k.$
For each $i=2,3,\ldots, k-1$, we merge the last vertex $t_i+2$ of $\tau(P_c^{(i)})$ and the first vertex $1$ of $\tau(P_c^{(i+1)})$ and relabel the vertices by $\{t_1 +2,t_1 +3,\ldots,p_i\}$ in increasing order. Let the resulting  noncrossing linked partition be denoted by $\nu (P_c)$. Then set
$$\varphi(P_i)=\{1,p_i+1\}\cup \varphi(P_c^{(1)})\cup\nu (P_c).$$
This operation is illustrated by Figure~\ref{Case3b}.

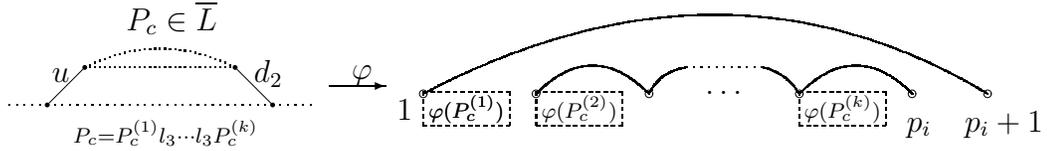
\begin{figure}[h]
\begin{center}
\setlength{\unitlength}{0.5cm}
\begin{picture}(23,4)(0,0)
\put(0,1){\circle*{.1}}\put(0.13,1.6){$u$}\put(0,1){\line(1,1){1}}\put(1,2){\circle*{.1}}\qbezier[30](1,2)(3,2)(5,2)
\qbezier[30](1,2)(3,3)(5,2)\put(5,2){\circle*{.1}}\put(5,2){\line(1,-1){1}}\put(6,1){\circle*{.1}}
\put(2.1,3){$P_c\in \overline{L}$}\put(5.5,1.6){$d_2$}\qbezier[40](-1,1)(3,1)(7,1)
\put(0.7,0){$\scriptstyle P_c=P_c^{(1)}l_3\cdots
l_3P_c^{(k)}$}

\put(7.5,1.5){\vector(1,0){1.5}}\put(8.1,1.7){$\varphi$}

\multiput(10,1.3)(3,0){3}{\circle{.2}}
\multiput(20,1.3)(3,0){2}{\circle{.2}}\put(25,1.3){\circle{.2}}
\put(9.3,0.6){$1$}
\put(22.85,0.3){$p_i$}\put(24.4,0.3){$p_i+1$}
\put(17.5,1.3){$\ldots$}

\qbezier(13,1.3)(14.5,2.8)(16,1.3)
\qbezier(16,1.3)(16.3,1.9)(17,2)\qbezier[10](17,2)(18,2)(19,2)
\qbezier(19,2)(19.6,1.9)(20,1.3)\qbezier(20,1.3)(21.5,2.8)(23,1.3)
\qbezier(10,1.3)(18,5.5)(25,1.3)

\put(10,0.45){\dashbox{.1}(2.3,.9){{$\scriptstyle
\varphi(P_c^{(1)})$}}}

\put(10,.45){\dashbox{.1}(2.3,.9){{$\scriptstyle
\varphi(P_c^{(1)})$}}}

\put(13,.45){\dashbox{.1}(2.3,.9){{$\scriptstyle
\varphi(P_c^{(2)})$}}}

\put(20,.45){\dashbox{.1}(2.3,.9){{$\scriptstyle
\varphi(P_c^{(k)})$}}}

\end{picture}
\caption{Case 3 (ii).}\label{Case3b}
\end{center}
\end{figure}
\end{itemize}

\end{enumerate}}

Finally, $\pi=\varphi(P)$ is constructed by merging  the last vertex of $\varphi(P_i)$ and the first vertex of $\varphi(P_{i+1})$, for $i=1,2,\ldots,k-1$, and relabeling the vertices by $\{1,2,\ldots,n+1\}$. It can be seen that $\pi$ is a noncrossing linked partition of $[n+1]$.

To show that $\varphi$ is a bijection, we give the inverse map of $\varphi$. Let $\pi\in NCL(n+1)$. We still work with the linear
representation of $\pi$.  First, we make use
of the outer arc decomposition of $\pi$. Here an outer arc is an arc in the linear representation of $\pi$ that is not covered  by any other arc.
To be more specific, the outer arc decomposition of $\pi$ is given by   $$\pi=(\pi_1,\pi_2,\ldots,\pi_m),$$ where each  $\pi_i$ is a noncrossing linked partition  of the set $\{s_i,s_i+1,\ldots,s_{i+1}-1,s_{i+1}\}$ with $s_1=1$, $s_i<s_{i+1}$ and $s_{m+1}=n+1$, such that $\pi_i=\{s_i\}\{s_{i+1}\}$, or $s_i$
and $s_{i+1}$ are contained in the same block of $\pi_i$. Next, consider $\varphi^{-1}(\pi_i)$, for $i=1,2,\ldots,m$. If $\pi_i=\{s_i,s_{i+1}\}$ (or $\{s_i\}\{s_{i+1}\}$), then set $\pi_i=l_1$ (or $l_2$). If $s_{i+1}\geq s_i+2$ and there are arcs forming a path from   the vertex $s_i$ to   the vertex $s_{i+1}-1$, then we deduce that $\varphi^{-1}(\pi_i)$ starts with an up step $u$ and ends  with a down step $d_1$. If there are no paths from the  vertex $s_i$ to the vertex $s_{i+1}-1$, then we have that $\varphi^{-1}(\pi_i)$ starts with an up step $u$ and ends with a down step $d_2$. Hence  $\varphi^{-1}(\pi_i)$ can be
reconstructed recursively according to the cases for the map $\varphi$. Finally, putting all the
pieces together, that is, setting
$$\varphi^{-1}(\pi)=\varphi^{-1}(\pi_1)\varphi^{-1}(\pi_2)\cdots
\varphi^{-1}(\pi_m),$$
we are led to a large
$(3,2)$-Motzkin path of length $n$. This completes the proof. \qed

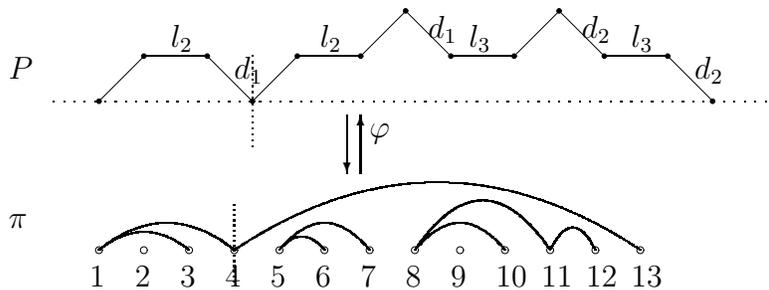
\begin{figure}[h]
\begin{center}
\setlength{\unitlength}{0.6cm}
\begin{picture}(14,6.4)(0,0)
\multiput(1,1)(1,0){13}{\circle{.15}}
\qbezier(1,1)(2,1.8)(3,1)\qbezier(1,1)(2.5,2.2)(4,1)\qbezier(4,1)(8.5,4)(13,1)
\qbezier(5,1)(5.5,1.6)(6,1)\qbezier(5,1)(6,2.2)(7,1)
\qbezier(8,1)(9,2.2)(10,1)\qbezier(8,1)(9.5,3.2)(11,1)\qbezier(11,1)(11.5,2)(12,1)
\put(.8,.2){$1$}\put(1.8,.2){$2$}\put(2.8,.2){$3$}\put(3.8,.2){$4$}
\qbezier[15](4,.3)(4,1.15)(4,2)\put(4.8,.2){$5$}
\put(5.8,.2){$6$}\put(6.8,.2){$7$}\put(7.8,.2){$8$}\put(8.8,.2){$9$}\put(9.8,.2){$10$}
\put(10.8,.2){$11$}\put(11.8,.2){$12$}\put(12.8,.2){$13$}\put(-1,1.5){$\pi$}

\put(6.5,4){\vector(0,-1){1.3}}\put(6.8,2.7){\vector(0,1){1.3}}
\put(7,3.5){$\varphi$}

\qbezier[70](0,4.3)(8,4.3)(16,4.3)\put(-1,4.8){$P$}
\put(1,4.3){\circle*{.1}}\put(1,4.3){\line(1,1){1}}\put(2,5.3){\circle*{.1}}
\put(2,5.3){\line(1,0){1.4}}\put(2.6,5.5){$l_2$}\put(3.4,5.3){\circle*{.1}}
\put(3.4,5.3){\line(1,-1){1}}\put(4.0,4.8){$d_1$}\qbezier[15](4.4,3.3)(4.4,4.3)(4.4,5.3)\put(4.4,4.3){\circle*{.1}}
\put(4.4,4.3){\line(1,1){1}}\put(5.4,5.3){\circle*{.1}}\put(5.4,5.3){\line(1,0){1.4}}
\put(5.9,5.5){$l_2$}\put(6.8,5.3){\circle*{.1}}\put(6.8,5.3){\line(1,1){1}}
\put(7.8,6.3){\line(1,-1){1}}\put(8.3,5.7){$d_1$}\put(7.8,6.3){\circle*{.1}}\put(8.8,5.3){\circle*{.1}}
\put(8.8,5.3){\line(1,0){1.4}}\put(9.2,5.5){$l_3$}\put(10.2,5.3){\circle*{.1}}
\put(10.2,5.3){\line(1,1){1}}\put(11.2,6.3){\circle*{.1}}\put(11.2,6.3){\line(1,-1){1}}
\put(12.2,5.3){\circle*{.1}}\put(11.7,5.7){$d_2$}\put(12.2,5.3){\line(1,0){1.4}}\put(12.8,5.5){$l_3$}
\put(13.6,5.3){\circle*{.1}}\put(13.6,5.3){\line(1,-1){1}}\put(14.2,4.8){$d_2$}\put(14.6,4.3){\circle*{.1}}
\end{picture}
\caption{Bijection $\varphi\colon L(12)\rightarrow
NCL(13)$.}\label{eg1}
\end{center}
\end{figure}

An example of the above bijection is given in Figure~\ref{eg1}.

The above bijection implies that the
large Schr\"{o}der number $S_{n}$ equals
the number $L_n$ of large $(3,2)$-Motzkin paths
of length $n$. On the other hand, there is a one-to-one correspondence between $(3,2)$-Motzkin paths
of length $n-1$ and little Schr\"{o}der paths of length $2n$. Therefore, the relation $S_{n}=2s_n$ can be rewritten as
\begin{equation}\label{lm}
 L_n=2m_{n-1},\end{equation}
that is, the number of
large $(3,2)$-Motzkin paths of length $n$ is twice the number of ordinary $(3,2)$-Motzkin paths of length $n-1$. Here we give a combinatorial interpretation of this fact.

Let $P$ be  a $(3,2)$-Motzkin path of length $n-1$.
If $P$ does not have any level step $l_3$ on the $x$-axis, then we get two large
$(3,2)$-Motzkin paths by adding a level
step $l_1$ or $l_2$ at the end of $P$.
Otherwise, remove the  first level step $l_3$  on the
$x$-axis in $P$,  and elevate the path after this $l_3$ level step. Concerning the elevated $(3,2)$-Motzkin path, there are two choices for the last down step.
It is easy to see that the above construction is reversible. Hence we obtain  (\ref{lm}).

\vskip 15pt

\noindent  {\bf Acknowledgments.}  This work was supported by  the 973
Project, the PCSIRT Project of the Ministry of Education,  and the National Science Foundation of China.

\end{document}